\documentclass[a4paper,10pt]{amsart}
\usepackage{latexsym}
\usepackage{amssymb}
\usepackage{mathrsfs}
\usepackage{amscd}
\usepackage[all]{xy}

\newtheorem{thm}{Theorem}[section]
\newtheorem{lem}[thm]{Lemma}
\newtheorem{prop}[thm]{Proposition}

\newtheorem*{Thm}{Theorem}

\theoremstyle{definition}
\newtheorem{defn}[thm]{Definition}
\newtheorem{ex}[thm]{Example}

\theoremstyle{remark}
\newtheorem{rem}[thm]{Remark}

\numberwithin{equation}{section}

\newcommand{\C}{{\mathbb C}}
\newcommand{\R}{{\mathbb R}}
\newcommand{\Z}{{\mathbb Z}}

\def\CA{\mathcal{A}}
\def\CB{\mathcal{B}}

\def\CT{\mathcal{T}}
\def\CU{\mathcal{U}}

\def\Fo{\mathfrak{o}}

\def\om{\omega}

\def\Aut{\operatorname{Aut}}
\def\Homeo{\operatorname{Homeo}}
\def\Isom{\operatorname{Isom}}
\def\id{\operatorname{id}}

\def\pr{\operatorname{pr}}

\newcommand{\abs}[1]{\lvert#1\rvert}

\begin{document}

\title[On Liftings of Local Torus Actions to Fiber bundles]{On Liftings of Local Torus Actions to Fiber Bundles}

\author[T. Yoshida]{Takahiko Yoshida}
\address{Graduate School of Mathematical Sciences, The University of Tokyo, \endgraf 
8-1 Komaba 3-chome, Meguro-ku, Tokyo, 153-8914, Japan}
\email{takahiko@ms.u-tokyo.ac.jp}
\thanks{The author is supported by Fujyukai Foundation and the 21st century COE program.}


\dedicatory{Dedicated to Professor Akio Hattori on his seventy seventh birthday.}


\begin{abstract}
In \cite{Y2} we introduced the notion of a local torus actions modeled on the standard representation (we call it a local torus action for simplicity), which is a generalization of a locally standard torus action. In this note we define a {\it lifting} of a local torus action to a principal torus bundle, and show that there is an obstruction class for the existence of liftings in the first cohomology of the fundamental group of the orbit space with coefficients in a certain module. 
\end{abstract}

\maketitle


\section{Introduction}
Let $G$ be a compact Lie group acting on a connected manifold $X$ and $\pi_P\colon P\to X$ a principal torus bundle on $X$. We denote the group of homeomorphisms of $X$ by $\Homeo (X)$ and identify the $G$-action on $X$ with the homomorphism $\phi \colon G\to \Homeo (X)$. We also denote the group of bundle isomorphisms of $P$ onto itself by $\Isom (P)$. Let $q \colon \Isom (P)\to \Homeo (X)$ be the obvious projection. We put $\Isom_G(P):=q^{-1}(\phi (G))$ and $\Aut (P):=\ker q$. Then there is an exact sequence of groups 
\[
1\to \Aut (P)\to \Isom_G(P)\xrightarrow{q} \phi (G) .
\]
A $G$-action $\tilde{\phi}\colon G\to \Isom_G(P)$ on $P$ as bundle isomorphisms which satisfies $q\circ \tilde{\phi}=\phi$ is called a {\it lifting} of $\phi$. It often comes into question whether $P$ admits a lifting of the $G$-action $\phi$ on $X$ or not. In \cite{St}, Stewart proved that $P$ admits a lifting if $G$ is simply connected and semi-simple. Hattori-Yoshida also proved in \cite{HY} that $P$ admits a lifting if and only if the Chern class of $P$ lies in the image of the equivariant cohomology group. For more details, consult \cite{HY}. 

In \cite{Y2}, as a generalization of locally standard torus actions and also an underlying structure of locally toric Lagrangian fibrations, we introduced the notion of a {\it local torus action modeled on the standard representation}  and proved the classification theorem for them. We also investigated their topology. The content of \cite{Y2} is a refinement of the work~\cite{Y1} and the talk by the author in International Conference on Toric Topology. In this note we discuss the lifting problem of local torus actions modeled on the standard representation in principal torus bundles. One of the motivation of this work is to generalize several works on the geometric quantization of symplectic toric manifolds to the case of locally toric Lagrangian fibrations. First we recall the definition of a local torus action modeled on the standard representation. Let $S^1$ be the unit circle in $\C$ and $T^n:=(S^1)^n$ the $n$-dimensional compact torus. $T^n$ acts on the $n$-dimensional complex vector space $\C^n$ by coordinatewise complex multiplication. This action is called the {\it standard representation of $T^n$}. 
\begin{defn}\label{localaction}
Let $X$ be a paracompact Hausdorff space. A {\it weakly standard $C^r$ $(0\le r\le \infty )$ atlas} of $X$ is an atlas $\{ (U_{\alpha}^X, \varphi^X_{\alpha})\}_{\alpha \in \CA}$ which satisfies the following properties 
\begin{enumerate}
\item for each $\alpha$, $\varphi^X_{\alpha}$ is a homeomorphism from $U_{\alpha}^X$ to an open set of $\C^n$ invariant under the standard representation of $T^n$,  
\item for each nonempty overlap $U_{\alpha \beta}^X:=U_{\alpha}^X\cap U_{\beta}^X$, 
\begin{enumerate}
\item $\varphi^X_{\alpha}(U_{\alpha \beta}^X)$ and $\varphi^X_{\beta}(U_{\alpha \beta}^X)$ are also invariant under the standard representation of $T^n$ and 
\item there exists an automorphism $\rho_{\alpha \beta}$ of $T^n$ as a Lie group such that the overlap map $\varphi^X_{\alpha \beta}:=\varphi^X_{\alpha}\circ (\varphi^X_{\beta})^{-1}$ is $\rho_{\alpha \beta}$-equivariant $C^r$ diffeomorphic with respect to the restrictions of the standard representation of $T^n$ to $\varphi^X_{\alpha}(U_{\alpha \beta}^X)$ and $\varphi^X_{\beta}(U_{\alpha \beta}^X)$. (The latter means that $\varphi^X_{\alpha \beta}(u\cdot z)=\rho_{\alpha \beta}(u)\cdot \varphi^X_{\alpha \beta}(z)$ for $u\in T^n$ and $z\in \varphi^X_{\beta}(U^X_{\alpha \beta})$.)
\end{enumerate}
\end{enumerate}
Two weakly standard $C^r$ atlases $\{ (U_{\alpha}^X, \varphi^X_{\alpha})\}_{\alpha \in \CA}$ and $\{ (V_{\beta}^X, \psi^X_{\beta})\}_{\beta \in \CB}$ of $X^{2n}$ are {\it equivalent} if on each nonempty overlap $U_{\alpha}^X\cap V_{\beta}^X$, there exists an automorphism $\rho$ of $T^n$ such that $\varphi^X_{\alpha}\circ (\psi^X_{\beta})^{-1}$ is $\rho$-equivariant $C^r$ diffeomorphic. We call an equivalence class of weakly standard $C^r$ atlases a {\it $C^r$ local $T^n$-action on $X^{2n}$ modeled on the standard representation} and denote it by $\CT$. 
\end{defn}
In the rest of this note, we called it a $C^r$ local $T^n$-action on $X$, or more simply, a local $T^n$-action on $X$ if there are no confusions. 

Next we define a lifting of a local torus action to a principal torus bundle. Let $(X, \CT)$ be a $2n$-dimensional manifold $X$ equipped with a local $T^n$-action $\CT$ and $\pi_P\colon P\to X$ a principal $T^k$-bundle on $X$. As we showed in \cite{Y2}, we can define the orbit space $B_X$ and the orbit map $\mu_X\colon X\to B_X$ for $(X, \CT)$ and the obstruction class in order that the local $T^n$-action on $X$ is induced by a global $T^n$-action lies in the first \v{C}ech cohomology $H^1(B_X;\Aut (T^n))$ of $B_X$ with coefficients in the group $\Aut (T^n )$ of group automorphisms of $T^n$. See Section~\ref{orbit-obst} for more details. Since $H^1(B_X;\Aut (T^n))$ is identified with the moduli space of representations of the fundamental group $\pi_1(B_X)$ of $B_X$ to $\Aut (T^n )$, by fixing a representative $\rho\colon \pi_1(B_X)\to \Aut (T^n)$ corresponding to the obstruction class, the fiber product $\pi^*X:=\{ (\tilde{b},x)\in \widetilde{B_X}\times X\colon \pi (\tilde{b})=\mu_X(x)\}$ of $\mu_X\colon X\to B_X$ and the universal covering $\pi\colon \widetilde{B_X}\to B_X$ of $B_X$ admits a global $T^n$-action $\phi_T\colon T^n\to \Homeo (\pi^*X)$. Moreover, by the construction, $\pi^*X$ is equipped with a natural $\pi_1(B_X)$-action $\phi_{\pi_1}\colon \pi_1(B_X)\to \Homeo (\pi^*X)$ and these actions form an action of the semidirect product $T^n\rtimes_{\rho}\pi_1(B_X)$ of $T^n$ and $\pi_1(B_X)$ with respect to $\rho$ on $\pi^*X$. For the explicit description of the action of the semidirect product, see Section~\ref{untwisting}. We put $G:=T^n\rtimes_{\rho}\pi_1(B_X)$. Let $\pi_{\widetilde{P}}\colon\widetilde{P}\to \pi^*X$ be the pullback of $\pi_P\colon P\to X$ to $\pi^*X$. $\widetilde{P}$ also admits a natural lifting $\tilde{\phi}_{\pi_1}\colon \pi_1(B_X)\to \Isom_G(\widetilde{P})$ of $\phi_{\pi_1}$. 
\begin{defn}\label{defn-lifting}
A {\it lifting} of the local $T^n$-action $\CT$ on $X$ to $P$ is a lifting $\tilde{\phi}_T\colon T^n\to \Isom_G(\widetilde{P})$ of $\phi_T$ such that 
\begin{equation}\label{semidirect-prod-condi}
\tilde{\phi}_T(\rho (a)(u))\circ \tilde{\phi}_{\pi_1}(a)=\tilde{\phi}_{\pi_1}(a)\circ \tilde{\phi}_T(u)
\end{equation}
for any $(u,a)\in G$. 
\end{defn}
In Proposition~\ref{lifting} we shall give an equivalent description of a lifting of $\CT$ to $P$ in terms of a sufficiently small weakly standard atlas $\{ (U_{\alpha}^X, \varphi^X_{\alpha})\}_{\alpha \in \CA}$ and the $T^n$-actions on $\varphi^X_{\alpha}(U_{\alpha}^X)$s. 

For the existence of liftings of $\phi_T$ which do not necessarily satisfy the condition \eqref{semidirect-prod-condi}, Hattori-Yoshida gave the necessary and sufficient condition in \cite{HY}. See also Theorem~\ref{Hattori-Yoshida}. We assume the existence of liftings of $\phi_T$ which do not necessarily satisfy the condition \eqref{semidirect-prod-condi}. Then the purpose of this note is to prove the following theorem. See Theorem~\ref{main-thm} for a precise statement. 
\begin{Thm}
There is an obstruction class in the first cohomology of the fundamental group $\pi_1(B_X)$ of $B_X$ with coefficients in a certain $\pi_1(B_X)$-module in order that $\pi_P\colon P\to X$ admits a lifting of the local $T^n$-action $\CT$ on $X$.  
\end{Thm}

This paper is organized as follows. In the next section, we recall the orbit space of a local torus action and the obstruction class in order that the local torus action is induced by a global torus action. We also give examples of a local torus action. In Section~\ref{untwisting} we explain that the fiber product of the orbit map of a local torus action and the universal covering of the orbit space admits a global torus action. We also describe the global torus action explicitly when a representation of $\pi_1(B_X)$ to $\Aut (T^n)$ corresponding to the obstruction class is given. Finally, in Section~\ref{lifting-problem}, we discuss the lifting problem of local torus actions in principal torus bundles and prove the above theorem. 

Throughout this paper we employ the vector notation in order to represent elements of $\C^n$, 
namely, $z=(z_1, \ldots ,z_n) \in \C^n$. The similar notation is also used for $T^n=(S^1)^n$, $\R^n$, etc. 

The author is thankful to the organizers of International Conference on Toric Topology for their hospitality. The author is also thankful for the helpful comments of the referee. 

\section{The orbit spaces and the obstruction to a global torus action}\label{orbit-obst}
Let $(X, \CT)$ be a $2n$-dimensional manifold $X$ equipped with a $C^r$ local $T^n$-action $\CT$. First we recall the orbit space $B_X$ of the local $T^n$-action $\CT$ on $X$. The orbit space $\C^n/T^n$ of the standard representation of $T^n$ is endowed with the natural stratification whose $k$-dimensional stratum consists of $k$-dimensional orbits. Let $\R^n_+$ be the standard $n$-dimensional positive cone 
\[
\R^n_+:=\{ \xi =(\xi_1, \ldots ,\xi_n)\in \R^n \colon \xi_i\ge 0\ i=1, \ldots ,n\} . 
\]
It also has a natural stratification with respect to the number of coordinates $\xi_i$ which are equal to zero. 
We define the map $\mu_{\C^n}\colon \C^n\to \R^n_+$ by 
\[
\mu_{\C^n}(z)=(\abs{z_1}^2, \ldots ,\abs{z_n}^2)
\]
for $z=(z_1, \ldots ,z_n)\in \C^n$. It is invariant under the standard representation of $T^n$ and induces a homeomorphism from $\C^n/T^n$ to $\R^n_+$ which preserves stratifications. 

Let $\{ (U_{\alpha}^X, \varphi^X_{\alpha})\}_{\alpha \in \CA}$ be a maximal weakly standard atlas of $X$ which belongs to $\CT$. We endow each quotient space $\varphi^X_{\alpha}(U_{\alpha}^X)/T^n$ with the quotient topology induced from the topology of $\varphi^X_{\alpha}(U_{\alpha}^X)$ by the natural projection $\pi_{\alpha} \colon \varphi^X_{\alpha}(U_{\alpha}^X)\to \varphi^X_{\alpha}(U_{\alpha}^X)/T^n$. By the property (2) for each overlap $U_{\alpha \beta}^X$, $\varphi^X_{\alpha \beta}$ induces a homeomorphism from $\varphi^X_{\beta}(U_{\alpha \beta}^X)/T^n$ to $\varphi^X_{\alpha}(U_{\alpha \beta}^X)/T^n$. We define two elements $b_{\alpha}\in \varphi^X_{\alpha}(U_{\alpha}^X)/T^n$ and $b_{\beta}\in \varphi^X_{\beta}(U_{\beta}^X)/T^n$ to be equivalent if $b_{\alpha}\in \varphi^X_{\alpha}(U_{\alpha \beta}^X)/T^n$, $b_{\beta} \in \varphi^X_{\beta}(U_{\alpha \beta}^X)/T^n$ and the map induced by $\varphi^X_{\alpha \beta}$ sends $b_{\beta}$ to $b_{\alpha}$. It is an equivalence relation on the disjoint union $\coprod_{\alpha}\left( \varphi^X_{\alpha}(U_{\alpha}^X)/T^n\right)$. We call the quotient space of $\coprod_{\alpha}\left( \varphi^X_{\alpha}(U_{\alpha}^X)/T^n\right)$ by the equivalence relation together with a quotient topology the {\it orbit space} of the local $T^n$-action $\CT$ on $X$ and denote it by $B_X$. It is easy to see that $B_X$ is a Hausdorff space and $\{ \varphi^X_{\alpha}(U_{\alpha}^X)/T^n\}_{\alpha \in \CA}$ is an open covering of $B_X$. By the construction of $B_X$, the map $\coprod_{\alpha}\pi_{\alpha} \circ \varphi^X_{\alpha}\colon \coprod_{\alpha}U_{\alpha}^X\to \coprod_{\alpha}\left( \varphi^X_{\alpha}(U_{\alpha}^X)/T^n\right)$ induces the map from $X$ to $B_X$. We call it the {\it orbit map} of the local $T^n$-action $\CT$ on $X$ and denote it by $\mu_X \colon X\to B_X$. Note that by the construction, it is a continuous open map. 
\begin{prop}\label{B}
$B_X$ is endowed with a structure of an $n$-dimensional topological manifold with corners. This means that $B_X$ has a system of coordinate neighborhoods modeled on open subsets of $\R^n_+$ so that overlap maps are homeomorphisms which preserve the natural stratifications induced from the one of $\R^n_+$. See \cite[Section 6]{Dav} for a topological manifold with corners. 
\end{prop}
\begin{proof}
We put $U_{\alpha}^B:=\varphi^X_{\alpha}(U_{\alpha}^X)/T^n$. The restriction of $\mu_{\C^n}$ to $\varphi^X_{\alpha}(U_{\alpha}^X)$ induces the homeomorphism from $U_{\alpha}^B$ to the open subset $\mu_{\C^n}(\varphi^X_{\alpha}(U_{\alpha}^X))$ of $\R^n_+$, which is denoted by $\varphi^B_{\alpha}$. By the construction, on each overlap $U_{\alpha \beta}^B:=U_{\alpha}^B\cap U_{\beta}^B$, the overlap map $\varphi^B_{\alpha \beta}:=\varphi^B_{\alpha}\circ (\varphi^B_{\beta})^{-1}\colon \mu_{\C^n}(\varphi^X_{\beta}(U_{\alpha \beta}^X))\to \mu_{\C^n}(\varphi^X_{\alpha}(U_{\alpha \beta}^X))$ preserves the natural stratifications of $\mu_{\C^n}(\varphi^X_{\alpha}(U_{\alpha \beta}^X))$ and $\mu_{\C^n}(\varphi^X_{\beta}(U_{\alpha \beta}^X))$. Thus, $\{ (U_{\alpha}^B,\varphi^B_{\alpha})\}_{\alpha \in \CA}$ is the desired atlas. 
\end{proof}
\begin{rem}\label{std}
The atlas $\{ (U_{\alpha}^B,\varphi^B_{\alpha})\}_{\alpha \in \CA}$ of $B_X$ constructed in the proof of Proposition~\ref{B} has following properties 
\begin{enumerate}
\item for each $\alpha$, $U_{\alpha}^X=\mu_X^{-1}(U_{\alpha}^B)$, 
$\varphi^X_{\alpha}(U_{\alpha}^X)=\mu_{\C^n}^{-1}(\varphi^B_{\alpha}(U_{\alpha}^B))$ 
and the following diagram commutes
\[
\xymatrix{
X\ar[d]^{\mu_X}\ar@{}[r]|{\supset}& \mu_X^{-1}(U_{\alpha}^B)\ar[r]^{\varphi^X_{\alpha}}\ar[d]^{\mu_X} & 
\mu_{\C^n}^{-1}(\varphi^B_{\alpha}(U_{\alpha}^B))\ar[d]^{\mu_{\C^n}}\ar@{}[r]|{\subset}& 
\C^n\ar[d]^{\mu_{\C^n}} \\
B_X\ar@{}[r]|\supset & U_{\alpha}^B\ar[r]^{\varphi^B_{\alpha}} & \varphi^B_{\alpha}(U_{\alpha}^B)\ar@{}[r]|\subset& 
\R^n_+,   
}
\]
\item the restriction of $\{ (U_{\alpha}^B,\varphi^B_{\alpha})\}_{\alpha \in \CA}$ to the interior 
$B_X\setminus \partial B_X$ of $B_X$ is a $C^r$ atlas of $B_X\setminus \partial B_X$. 
\end{enumerate}
\end{rem}

We give some examples. 
\begin{ex}[Locally standard torus actions]\label{locallystandardtorusaction}
Let $T^n$ act smoothly on a $2n$-dimensional smooth manifold $X$. A {\it standard coordinate neighborhood} of $X$ consists of a triple $(U, \rho , \varphi)$, where $U$ is a $T^n$-invariant open set of $X$, $\rho$ is an automorphism of $T^n$, and $\varphi$ is a $\rho$-equivariant diffeomorphism from $U$ to some $T^n$-invariant open subset in $\C^n$. The action of $T^n$ on $X$ is said to be {\it locally standard} if every point in $X$ lies in some standard coordinate neighborhood. See \cite{DJ, BP} for more details. (A typical example of locally standard torus actions is a nonsingular toric variety.) The atlas which consists of standard coordinate neighborhoods is weakly standard. Hence, a locally standard $T^n$-action induces a local $T^n$-action on $X$. 
\end{ex}

Note that not all local torus actions are induced by locally standard torus actions. For example, a nontrivial $T^n$-bundle on an $n$-dimensional closed manifold whose structure group is $\Aut (T^n)$ is equipped with a local $T^n$-action which is not induced by any locally standard $T^n$-action. In general, for any $C^r$ local $T^n$-action $\CT$ on a $2n$-dimensional manifold $X$, we take a weakly standard atlas $\{ (U_{\alpha}^X, \varphi^X_{\alpha})\}_{\alpha \in \CA}$ belonging to $\CT$. We put $\CU:=\{ U_{\alpha}^B\}_{\alpha \in \CA}$. It is easy to see that the automorphisms $\rho_{\alpha \beta}$ of $T^n$ in the property (2) of Definition~\ref{localaction} form a \v{C}ech one-cocycle $\{ \rho_{\alpha \beta}\}$ on $\CU$ with values in $\Aut (T^n)$. 
\begin{prop}\label{ob1}
A $C^r$ local $T^n$-action on $X$ is induced by some $C^r$ locally standard $T^n$-action if and only if 
$\{ \rho_{\alpha \beta}\}$ and the trivial \v{C}ech one-cocycle are of the same equivalence class in the first \v{C}ech cohomology set $H^1(B_X;\Aut (T^n))$, where the trivial \v{C}ech one-cocycle is the one whose values on all open set are equal to the identity map of $T^n$.  
\end{prop}
For the proof, see \cite{Y2}. 
\begin{rem}\label{rem-rep}
It is well known that there is a one-to-one correspondence between $H^1(B_X;\Aut (T^n))$ and the moduli space of representations of $\pi_1(B_X)$ to $\Aut (T^n)$. 
\end{rem}

\begin{ex}\label{ex-cylinder}
We can construct an example of local torus actions which does not come from any locally standard torus fibrations in the following way. Let $\overline{X}$ be the quotient space of the space 
\[
\{ (\xi ,u ,z)\in \R^2\times T^2\times \C^2 \colon \xi_2=\abs{z_1}^2, \xi_2+\abs{z_2}^2=1\}
\]
by the $T^2$-action defined by 
\[
v\cdot (\xi ,u ,z):=\left( \xi , (u_1, u_2v_1v_2^{-1}), v^{-1}\cdot z\right) . 
\]
There is a $T^2$-action on $\overline{X}$ which is induced by the multiplication to the second factor of $\R^2\times T^2\times \C^2$. The orbit space of this action is naturally identified with $\overline{B}:=\R\times [0,1]$ and the orbit map is induced by the first projection of $\R^2\times T^2\times \C^2$. Define the right action of $\Z$ on $\overline{X}$ by
\[
[\xi , u ,z]\cdot n =[(\xi_1-n(\xi_2+1), \xi_2), \rho^{-n}(u), z]
\]
for $(\xi , \theta ,z)\in \overline{X}$ and $n\in \Z$, where $\rho^n(u):=(u_1,u_1^{-n}u_2)$. It is easy to see that the $\Z$-action on $\overline{X}$ is well-defined, and the $\Z$-action and the $T^2$-action on $\overline{X}$ satisfies the following condition
\[
(v\cdot [\xi , u ,z])\cdot n=\rho^{-n}(v)\cdot ([\xi , u ,z]\cdot n)
\]
for $v\in T^2$ and $n\in \Z$. Moreover, the $\Z$-action on $\overline{X}$ descends to the action of $\Z$ on $\overline{B}$ which is defined by 
\begin{equation}\label{deck-trans}
\xi \cdot n=(\xi_1-n(\xi_2+1), \xi_2). 
\end{equation}
We denote by $X$, $B$ the quotient spaces $\overline{X}/\Z$, $\overline{B}/\Z$, respectively. We also denote by $\mu \colon X\to B$ the map induced by the orbit map of the $T^n$-action on $\overline{X}$. Then, it is easy to see that $X$ is equipped with a local $T^2$-action whose orbit space is $B$ and the orbit map is $\mu$.  
\end{ex} 
\begin{ex}[Locally toric Lagrangian fibrations~\cite{Ham}]
One of the important examples is a locally toric Lagrangian fibration. Let $(X,\om)$ be a $2n$-dimensional symplectic manifold and $B$ an $n$-dimensional smooth manifold with corners. A map $\mu \colon X\to B$ is called a {\it locally toric Lagrangian fibration} if for each point $b\in B$, there exists a coordinate neighborhood $(U,\varphi^B)$ of $b$ modeled on $\R^n_+$ and there exists a symplectomorphism $\varphi^X\colon (\mu^{-1}(U), \om )\to (\mu_{\C^n}^{-1}(\varphi^B(U)), \om_{\C^n})$ such that $\mu_{\C^n}\circ \varphi^X=\varphi^B\circ \mu$. A locally toric Lagrangian fibration is endowed with a smooth local $T^n$-action. See \cite{Y2}, for more details. 
\end{ex}

\section{Untwisting local torus actions}\label{untwisting}
Let $(X,\CT )$ be a $2n$-dimensional manifold $X$ equipped with a local $T^n$-action $\CT$, $\{ (U_{\alpha}^X, \varphi^X_{\alpha})\}_{\alpha \in \CA}$ a weakly standard atlas of $X$ belonging to $\CT$, and $\{ (U_{\alpha}^B, \varphi^B_{\alpha})\}_{\alpha \in \CA}$ the atlas of $B_X$ induced by $\{ (U_{\alpha}^X, \varphi^X_{\alpha})\}_{\alpha \in \CA}$ which satisfies the properties of Remark~\ref{std}. Let $\pi\colon \widetilde{B_X}\to B_X$ be the universal covering of $B_X$. It is obvious that the fiber product $\pi^*X$ of $\mu_X\colon X\to B_X$ and $\pi\colon \widetilde{B_X}\to B_X$ admits a local $T^n$-action whose orbit space is $\widetilde{B_X}$. For more details, see \cite[Example 3.14]{Y2}. Since $\widetilde{B_X}$ is simply connected, by Proposition~\ref{ob1} and Remark~\ref{rem-rep}, this local $T^n$-action comes from a locally standard $T^n$-action. Moreover, by the construction, $\pi^*X$ is equipped with the natural action of $\pi_1(B_X)$ and these two actions form an action of the semidirect product of $T^n$ and $\pi_1(B_X)$. The aim of this section is to give the explicit description of this action . 

By replacing $\{ U_{\alpha}^B\}$ by its refinement if necessary, we may assume that for each $\alpha$, there exists a local trivialization $\varphi^{\tilde{B}}_{\alpha}\colon \pi^{-1}(U_{\alpha}^B)\to U_{\alpha}^B\times \pi_1(B_X)$ of $\pi\colon \widetilde{B_X}\to B_X$ as a principal $\pi_1(B_X)$-bundle. On each nonempty overlap $U_{\alpha \beta}^B$, we denote the transition function with respect to these local trivializations by $a_{\alpha \beta}$. Note that $a_{\alpha \beta}$ is locally constant since $\pi_1(B_X)$ is discrete. 

As we described before, the automorphisms $\rho_{\alpha \beta}$ of $T^n$ in the property (2) of Definition~\ref{localaction} form a \v{C}ech cohomology class $[\{ \rho_{\alpha \beta}\}]\in H^1(B_X;\Aut (T^n))$. We take a representative $\rho \colon \pi_1(B_X)\to \Aut (T^n)$ of the equivalence class of representations corresponding to $[\{ \rho_{\alpha \beta}\}]$. Note that $\rho$ is unique up to the conjugation of $\Aut (T^n)$. Then, for each $\alpha$ there exists an automorphism $\rho_{\alpha}\in \Aut (T^n)$ such that $\rho_{\alpha \beta}=\rho_{\alpha}\circ \rho (a_{\alpha \beta})\circ \rho_{\beta}^{-1}$ on each nonempty overlap $U_{\alpha \beta}^B$. Let $T^n\rtimes_{\rho} \pi_1(B_X)$ be the semidirect product of $T^n$ and $\pi_1(B_X)$ with respect to $\rho$, namely, $T^n\rtimes_{\rho} \pi_1(B_X)$ is the Cartesian product of $T^n$ and $\pi_1(B_X)$ as a set with the product 
\[
(u_1,a_1)(u_2,a_2):=(u_1\rho (a_1)(u_2), a_1a_2) . 
\]
Let $(u,a)$ be an element of $T^n\rtimes_{\rho} \pi_1(B_X)$ and $(\tilde{b}, x)$ an element of $\pi^*X$. Suppose that $\tilde{b}$ lies in $\pi^{-1}(U_{\alpha}^B)$ and $\varphi^{\tilde{B}}_{\alpha}(\tilde{b})=(\pi(\tilde{b}), a_{\alpha})$. We put 
\begin{equation}\label{prod-semidirect} 
(u,a)\cdot_{\rho} (\tilde{b},x):=\left( \tilde{b}\cdot a^{-1}, (\varphi^X_{\alpha})^{-1}(\rho_{\alpha}\circ \rho(a_{\alpha}a^{-1})(u)\cdot \varphi^X_{\alpha}(x))\right) . 
\end{equation}
\begin{lem}\label{phi_G}
\eqref{prod-semidirect} does not depend on the choice of local trivializations and defines the action of $T^n\rtimes_{\rho} \pi_1(B_X)$ on $\pi^*X$. 
\end{lem}
\begin{proof} 
Suppose that $\tilde{b}$ also lies in $\pi^{-1}(U_{\beta}^B)$ for another $U_{\beta}^B$ and $\varphi^{\tilde{B}}_{\beta}(\tilde{b})=(\pi(\tilde{b}), a_{\beta})$. By using $a_{\alpha}=a_{\alpha \beta}a_{\beta}$ and $\rho_{\alpha \beta}=\rho_{\alpha}\circ \rho (a_{\alpha \beta})\circ \rho_{\beta}^{-1}$, 
\[
\begin{split}
&(\varphi^X_{\alpha})^{-1}(\rho_{\alpha}\circ \rho(a_{\alpha}a^{-1})(u)\cdot \varphi^X_{\alpha}(x))\\
&=(\varphi^X_{\beta})^{-1}\circ (\varphi^X_{\alpha \beta})^{-1}\left( \rho_{\alpha}\circ \rho(a_{\alpha \beta}a_{\beta}a^{-1})(u)\cdot \varphi^X_{\alpha \beta}\circ \varphi^X_{\beta}(x)\right) \\
&=(\varphi^X_{\beta})^{-1}\left( \rho_{\alpha \beta}^{-1}\circ \rho_{\alpha}\circ \rho(a_{\alpha \beta})\circ \rho(a_{\beta}a^{-1})(u)\cdot \varphi^X_{\beta}(x)\right) \\
&=(\varphi^X_{\beta})^{-1}(\rho_{\beta}\circ \rho(a_{\beta}a^{-1})(u)\cdot \varphi^X_{\beta}(x)). 
\end{split}
\]
This implies that \eqref{prod-semidirect} does not depend on the choice of local trivializations. Next, we check that \eqref{prod-semidirect} defines an action. For $(u_1, a_1)$, $(u_2, a_2)\in T^n\rtimes_{\rho} \pi_1(B_X)$ and $(\tilde{b},x)\in \pi^*X$ satisfying $\tilde{b}\in \pi^{-1}(U_{\alpha}^B)$ and $\varphi^{\tilde{B}}_{\alpha}(\tilde{b})=(\pi(\tilde{b}), a_{\alpha})$, 
\[
\begin{split}
&(u_1, a_1)\cdot_{\rho} \left( (u_2, a_2)\cdot_{\rho} (\tilde{b},x)\right) \\
&=(u_1, a_1)\cdot_{\rho} \left( \tilde{b}\cdot a_2^{-1}, (\varphi^X_{\alpha})^{-1}(\rho_{\alpha}\circ \rho(a_{\alpha}a_2^{-1})(u_2)\cdot \varphi^X_{\alpha}(x))\right) \\
&=\left( \tilde{b}\cdot a_2^{-1}a_1^{-1}, (\varphi^X_{\alpha})^{-1}(\rho_{\alpha}\circ \rho(a_{\alpha}a_2^{-1}a_1^{-1})(u_1)\rho_{\alpha}\circ \rho(a_{\alpha}a_2^{-1})(u_2)\cdot \varphi^X_{\alpha}(x))\right) \\ 
&=\left( \tilde{b}\cdot (a_1a_2)^{-1}, (\varphi^X_{\alpha})^{-1}(\rho_{\alpha}\circ \rho(a_{\alpha}(a_1a_2)^{-1})(u_1\rho(a_1)(u_2))\cdot \varphi^X_{\alpha}(x))\right) \\ 
&= (u_1\rho(a_1)(u_2), a_1a_2)\cdot_{\rho} (\tilde{b}, x). 
\end{split}
\]
This proves the lemma. 
\end{proof}
Note that the orbit space $(\pi^*X)/T^n\rtimes_{\rho} \pi_1(B_X)$ is naturally identified with $B_X$. 

Let $\rho'\colon \pi_1(B_X)\to \Aut (T^n)$ be another representative of the equivalence class of representations corresponding to $[\{ \rho_{\alpha \beta}\}]$. Then, there exists an automorphism $f\in \Aut (T^n)$ such that $\rho'=f\circ \rho \circ f^{-1}$. $f$ defines the group isomorphism $\overline{f}\colon T^n\rtimes_{\rho} \pi_1(B_X)\to T^n\rtimes_{\rho'} \pi_1(B_X)$ by $\overline{f}(u,a):=(f(u),a)$. Then we can check the following proposition. 
\begin{prop}
For any $(u,a)\in T^n\rtimes_{\rho} \pi_1(B_X)$ and $(\tilde{b},x)\in \pi^*X$, 
\[
(u,a)\cdot_{\rho}(\tilde{b},x)=\overline{f}(u,a)\cdot_{\rho'}(\tilde{b},x). 
\]
\end{prop}

\section{An obstruction theory}\label{lifting-problem}
In this section we give an obstruction class for the existence of a lifting when a local torus action on a manifold and a principal torus bundle on it are given. We assume that manifolds, maps, and local $T^n$-actions are of class $C^0$ unless otherwise stated. 

Let $(X,\CT )$ be a $2n$-dimensional manifold $X$ equipped with a local $T^n$-action $\CT$ and $\pi\colon \widetilde{B_X}\to B_X$ the universal covering of $B_X$. We fix a representative $\rho \colon \pi_1(B_X)\to \Aut (T^n)$ of the equivalence class of representations corresponding to $[\{ \rho_{\alpha \beta}\}]$ and denote $T^n\rtimes_{\rho} \pi_1(B_X)$ by $G$. Then, $G$ acts on $\pi^*X$ by \eqref{prod-semidirect}. In the rest of this paper we omit the subscript $_{\rho}$ of the $G$-action on $\pi^*X$ in \eqref{prod-semidirect}. 

Suppose that $\pi_P\colon P\to X$ is a principal $T^k$-bundle on $X$. Then the pullback $\pi_{\widetilde{P}}\colon \widetilde{P}\to \pi^*X$ of $\pi_P\colon P\to X$ to $\pi^*X$ can be written by 
\[
\widetilde{P}=\{ (\tilde{b},q)\in \widetilde{B_X}\times P\colon \pi(\tilde{b})=\mu_X\circ \pi_P(q)\}. 
\]
$\widetilde{P}$ admits a natural lifting $\tilde{\phi}_{\pi_1}\colon \pi_1(B_X)\to \Isom_G(\widetilde{P})$ of the natural $\pi_1(B_X)$-action $\phi_{\pi_1}\colon \pi_1(B_X)\to \Homeo (\pi^*X)$ on $\pi^*X$ which is written by 
\[
\tilde{\phi}_{\pi_1}(a)(\tilde{p}):=(\tilde{b}\cdot a^{-1}, p)
\]
for $\tilde{p}=(\tilde{b},p)\in \widetilde{P}$ and $a\in \pi_1(B_X)$. We also represent the $T^n$-action on $\pi^*X$ as the subgroup of $G$ by $\phi_T\colon T^n\to \Homeo (\pi^*X)$. 
\begin{prop}\label{lifting}
Let $\{ (U_{\alpha}^X, \varphi^X_{\alpha})\}_{\alpha \in \CA}$ be a weakly standard atlas of $X$ belonging to $\CT$ and $\{ (U_{\alpha}^B, \varphi^B_{\alpha})\}_{\alpha \in \CA}$ the atlas of $B_X$ induced by $\{ (U_{\alpha}^X, \varphi^X_{\alpha})\}_{\alpha \in \CA}$ which satisfies the properties of Remark~\textup{\ref{std}}. We assume that for each $\alpha$, there exists a local trivialization $\varphi^{\tilde{B}}_{\alpha}\colon \pi^{-1}(U_{\alpha}^B)\to U_{\alpha}^B\times \pi_1(B_X)$ of $\pi\colon \widetilde{B_X}\to B_X$. Then $P$ admits a lifting of the local $T^n$-action $\CT$ on $X$ if and only if there exists a family $\{ (P_{\alpha}, \tilde{\phi}_{\alpha}, \varphi^P_{\alpha})\}_{\alpha \in \CA}$ of triples, where for each $\alpha$, $(P_{\alpha}, \tilde{\phi}_{\alpha}, \varphi^P_{\alpha})$ consists of a principal $T^k$-bundle $\pi_{P_{\alpha}}\colon P_{\alpha}\to \varphi^X_{\alpha}(U_{\alpha}^X)$ equipped with a lifting $\tilde{\phi}_{\alpha}\colon T^n\to \Isom_{T^n}(P_{\alpha})$ of the $T^n$-action $\phi_{\alpha}\colon T^n\to \Homeo (\varphi^X_{\alpha}(U_{\alpha}^X))$ on $\varphi^X_{\alpha}(U_{\alpha}^X)$ obtained as the restriction of  the standard representation of $T^n$ and a bundle isomorphism $\varphi^P_{\alpha}\colon P|_{U_{\alpha}^X}\to P_{\alpha}$ which covers $\varphi^X_{\alpha}$ such that on each nonempty overlap $U_{\alpha \beta}^X$, $\varphi^P_{\alpha \beta}:=\varphi^P_{\alpha}\circ (\varphi^P_{\beta})^{-1}\colon P_{\beta}|_{\varphi^X_{\beta}(U_{\alpha \beta}^X)}\to P_{\alpha}|_{\varphi^X_{\alpha}(U_{\alpha \beta}^X)}$ is $\rho_{\alpha \beta}$-equivariant with respect to the liftings. 
\end{prop}
\begin{proof}
If there exists such a family $\{ (P_{\alpha}, \tilde{\phi}_{\alpha}, \varphi^P_{\alpha})\}_{\alpha \in \CA}$, then we can define a lifting of $\phi_T$ as follows. For $u\in T^n$ and $(\tilde{b},p)\in \widetilde{P}$ where $\tilde{b}\in \pi^{-1}(U_{\alpha}^B)$ and $\varphi^{\tilde{B}}_{\alpha}(\tilde{b})=(\pi(\tilde{b}), a_{\alpha})$, we define a lifting $\tilde{\phi}_T\colon T^n\to \Isom_G(\widetilde{P})$ of $\phi_T$ by 
\[
\tilde{\phi}_T(u)(\tilde{b},p):=\left( \tilde{b}, (\varphi^P_{\alpha})^{-1}(\tilde{\phi}_{\alpha}(\rho_{\alpha}\circ \rho(a_{\alpha})(u))(\varphi^P_{\alpha}(p)))\right) .
\]
We can show that it is well-defined by the same way as in the proof of Lemma~\ref{phi_G}. It is also easy to check that $\tilde{\phi}_T$ satisfies \eqref{semidirect-prod-condi}. 

Conversely, suppose that $\widetilde{P}$ admits a lifting $\tilde{\phi}_T$ of $\phi_T$ which satisfies \eqref{semidirect-prod-condi}. We put $P_{\alpha}:=P|_{U^X_{\alpha}}$, $\pi_{\alpha}:=\varphi^X_{\alpha}\circ \pi_P\colon P_{\alpha}\to \varphi^X_{\alpha}(U^X_{\alpha})$, and $\varphi^P_{\alpha}:=\id$ for each $\alpha$. For any point $p\in P_{\alpha}$, $\left( (\varphi^{\tilde{B}}_{\alpha})^{-1}(\mu_X\circ \pi_P(p),e),p\right)$ is an element of $\pi^*P$, where $e$ is the unit element of $\pi_1(B_X)$. Then, the following equation defines a lifting $\tilde{\phi}_{\alpha}\colon T^n\to \Isom_T(P_{\alpha})$ of $\phi_{\alpha}\colon T^n\to \Homeo (\varphi^X_{\alpha}(U_{\alpha}^X))$
\[
\tilde{\phi}_T(\rho_{\alpha}^{-1}(u))\left( (\varphi^{\tilde{B}}_{\alpha})^{-1}(\mu_X\circ \pi_P(p),e),p\right) =\left( (\varphi^{\tilde{B}}_{\alpha})^{-1}(\mu_X\circ \pi_P(p),e),\tilde{\phi}_{\alpha}(u)(p)\right) . 
\]
By using \eqref{semidirect-prod-condi} and $\rho_{\alpha \beta}=\rho_{\alpha}\circ \rho (a_{\alpha \beta})\circ \rho_{\beta}^{-1}$ on $U_{\alpha \beta}^B$, for any $p\in P|_{U_{\alpha \beta}^X}$, 
\begin{equation*}
\begin{split}
&\left( (\varphi^{\tilde{B}}_{\beta})^{-1}(\mu_X\circ \pi_P(p),e),\tilde{\phi}_{\beta}(u)(p)\right) \\
&=\tilde{\phi}_T(\rho_{\beta}^{-1}(u))\left( (\varphi^{\tilde{B}}_{\beta})^{-1}(\mu_X\circ \pi_P(p),e),p\right) \\
&=\tilde{\phi}_T(\rho_{\beta}^{-1}(u))\left( (\varphi^{\tilde{B}}_{\alpha})^{-1}\circ \varphi^{\tilde{B}}_{\alpha \beta}(\mu_X\circ \pi_P(p),e),p\right) \\
&=\tilde{\phi}_T(\rho_{\beta}^{-1}(u))\left( (\varphi^{\tilde{B}}_{\alpha})^{-1}(\mu_X\circ \pi_P(p),a_{\alpha \beta}),p\right) \\
&=\tilde{\phi}_T(\rho_{\beta}^{-1}(u))\circ \tilde{\phi}_{\pi_1}(a_{\alpha \beta}^{-1})\left( (\varphi^{\tilde{B}}_{\alpha})^{-1}(\mu_X\circ \pi_P(p),e),p\right) \\
&=\tilde{\phi}_{\pi_1}(a_{\alpha \beta}^{-1})\circ \tilde{\phi}_T(\rho(a_{\alpha \beta})\circ \rho_{\beta}^{-1}(u))\left( (\varphi^{\tilde{B}}_{\alpha})^{-1}(\mu_X\circ \pi_P(p),e),p\right)\\
&=\tilde{\phi}_{\pi_1}(a_{\alpha \beta}^{-1})\circ \tilde{\phi}_T(\rho_{\alpha}^{-1}\circ \rho_{\alpha \beta}(u))\left( (\varphi^{\tilde{B}}_{\alpha})^{-1}(\mu_X\circ \pi_P(p),e),p\right)\\
&=\tilde{\phi}_{\pi_1}(a_{\alpha \beta}^{-1})\left( (\varphi^{\tilde{B}}_{\alpha})^{-1}(\mu_X\circ \pi_P(p),e),\tilde{\phi}_{\alpha}(\rho_{\alpha \beta}(u))(p)\right)\\
&=\left( (\varphi^{\tilde{B}}_{\alpha})^{-1}(\mu_X\circ \pi_P(p),a_{\alpha \beta}),\tilde{\phi}_{\alpha}(\rho_{\alpha \beta}(u))(p)\right)\\
&=\left( (\varphi^{\tilde{B}}_{\beta})^{-1}(\mu_X\circ \pi_P(p),e),\tilde{\phi}_{\alpha}(\rho_{\alpha \beta}(u))(p)\right) . 
\end{split}
\]
This implies that $\varphi^P_{\alpha \beta}$ is $\rho_{\alpha \beta}$-equivariant, namely, 
\[
\varphi^P_{\alpha \beta}\circ \tilde{\phi}_{\beta}(u)=\tilde{\phi}_{\alpha}(\rho_{\alpha \beta}(u))\circ \varphi^P_{\alpha \beta}. 
\]
This proves the proposition. 
\end{proof}

\begin{ex}
We consider Example~\ref{ex-cylinder}. In this case, the universal covering of $B$ is $\overline{B}$ and the pullback of $X$ to $\overline{B}$ is naturally identified with $\overline{X}$. The fundamental group $\pi_1(B)$ is isomorphic to $\Z$ and the deck transformation is given by \eqref{deck-trans}. Under the identification $\pi_1(B)\cong \Z$, the representation of $\pi_1(B)$ to $\Aut (T^2)$ is given by assigning the element $\rho^n(v)\in T^2$ to $n\in \pi_1(B)$ and $v\in T^2$, which is also denoted by $\rho\colon \pi_1(B)\to \Aut (T^2)$. Then, the $T^2$-action and the right $\Z$-action on $\overline{X}$ form an action $\phi \colon G\to \Homeo (\overline{X})$ of the semidirect product $G=T^2\rtimes_{\rho}\pi_1(B)$ which is written by 
\[
\phi(v,n)([\xi ,u , z])=v\cdot ([\xi ,u ,z]\cdot (-n))=[(\xi_1+n(\xi_2+1), \xi_2), \rho^n(u)v, z]
\]
for $(v,n)\in G$ and $[\xi ,u, z]\in \overline{X}$. 

Let $\pi_{\widetilde{P}}\colon \widetilde{P}\to \overline{X}$ be a principal circle bundle which is defined by the quotient space
\[
\widetilde{P}:=\{ (\xi ,u ,z)\in \R^2\times T^2\times \C^2 \colon \xi_2=\abs{z_1}^2, \xi_2+\abs{z_2}^2=1\} \times_{T^2} S^1
\]
by the $T^2$-action 
\[
v\cdot (\xi ,u ,z, t):=\left( \xi , (u_1, u_2v_1v_2^{-1}), v^{-1}\cdot z, v_2t\right) . 
\]
We fix an element $s\in S^1$ and define the lifting $\tilde{\phi} \colon G\to \Isom_G (\widetilde{P})$ of $\phi$ to $\widetilde{P}$ by
\[
\tilde{\phi}(v,n)([\xi ,u ,z, t]):=[(\xi_1+n(\xi_2+1), \xi_2), \rho^n(u)v, z, v_1s^nt]
\]
for $(v,n)\in G$ and $[\xi ,u ,z, t]\in \widetilde{P}$. It is easy to see that $\tilde{\phi}$ is well-defined. In particular, the composition of the natural homomorphism $s\colon \pi_1(B)\to G$ defined by $s(n):=(1,n)$ and $\tilde{\phi}$ defines a $\pi_1(B)$-action on $\widetilde{P}$. We denote by $P$ the quotient space of this $\pi_1(B)$-action. By the construction, $\pi_{\widetilde{P}}\colon \widetilde{P}\to \overline{X}$ descends to the principal circle bundle $\pi_P\colon P\to X$ which admits a lifting of the local $T^2$-action on $X$ defined in Example~\ref{ex-cylinder}. 
\end{ex}

In the rest of this note, let us investigate when $P$ admits a lifting of the local $T^n$-action $\CT$ on $X$. For the existence of a lifting of $\phi_T$, which does not necessarily satisfy the condition \eqref{semidirect-prod-condi}, Hattori-Yoshida gave the necessary and sufficient condition in \cite{HY} which is described as follows. Let $(\pi^*X)_{T^n}$ be the Borel construction $ET^n\times_{\phi_T}\pi^*X$ with respect to the $T^n$-action $\phi_T$ on $\pi^*X$. We denote by $r\colon ET^n\times \pi^*X\to ET^n\times_{\phi_T}\pi^*X$ and $\pr_2\colon ET^n\times \pi^*X\to \pi^*X$ the natural projections. Since $ET^n$ is contractible, the induced homomorphism $\pr_2^*\colon H^2(\pi^*X; \Z^n)\to H^2(ET^n\times \pi^*X; \Z^n)$ is an isomorphism. 
\begin{thm}[\cite{HY}]\label{Hattori-Yoshida}
$\pi_{\widetilde{P}}\colon \widetilde{P}\to \pi^*X$ admits a lifting of the $T^n$-action $\phi_T$ on $\pi^*X$ if and only if the Chern class $c_1(\widetilde{P})$ lies in the image of ${\pr_2^*}^{-1}\circ r^*\colon H^2(ET^n\times_{\phi_T}\pi^*X; \Z^n)\to H^2(\pi^*X;\Z^n)$. 
\end{thm}

Next, we assume that $\pi_{\widetilde{P}}\colon \widetilde{P}\to \pi^*X$ satisfies the condition in Theorem~\ref{Hattori-Yoshida} and investigate when $\widetilde{P}$ admits a lifting of $\phi_T$ which satisfies \eqref{semidirect-prod-condi}. We fix a lifting $\tilde{\phi}_T\colon T^n\to \Isom_G(\widetilde{P})$ of $\phi_T$ to $\widetilde{P}$. Then $T^n$ acts on $\Aut (\widetilde{P})$ by the conjugation $\tilde{\phi}_T^{-1}(u)\circ f\circ \tilde{\phi}_T(u)$ for $u\in T^n$ and $f\in \Aut (\widetilde{P})$. Note that the group $\Aut (\widetilde{P})$ is canonically isomorphic to the group $C(\pi^*X,T^k)$ of all continuous maps from $\pi^*X$ to $T^k$. The isomorphism is given by assigning to each $f\in \Aut (\widetilde{P})$ the map $\tau \in C(\pi^*X,T^k)$ which is determined uniquely by the equation 
\[
f(\tilde{p})=\tilde{p}\cdot \tau (\pi_{\widetilde{P}}(\tilde{p})).
\]
For more details, consult \cite{HY}. Through the isomorphism from $\Aut (\widetilde{P})$ to $C(\pi^*X, T^k)$, $T^n$ also acts on $C(\pi^*X, T^k)$ by $\tau^u(\tilde{x}):=\tau (\phi_T(u)(\tilde{x}))$. In general, if $K$ is a group and $M$ is a topological right $K$-module, then we can define the cochain complex $C^*(K;M)$ of continuous cochains of $K$ with values in $M$ as follows. That is, $C^q(K; M)$ is the abelian group of all continuous maps from the $q$th Cartesian product $K^q$ of $K$ into $M$ and the coboundary $\delta \colon C^q(K; M)\to C^{q+1}(K; M)$ is defined by 
\[
\begin{split}
\delta \sigma (u_1, \ldots ,u_{q+1}):=&\sigma (u_2, \ldots ,u_{q+1})+
\sum_{i=1}^q(-1)^i\sigma (u_1,\ldots ,u_iu_{i+1}, \ldots ,u_{q+1})\\
&+(-1)^{q+1}\sigma (u_1,\ldots ,u_q)\cdot {u_{q+1}}. 
\end{split}
\]
Let $C^*(T^n; C(\pi^*X, T^k))$ be the cochain complex with values in the right $T^n$-module $C(\pi^*X, T^k)$. $\pi_1(B_X)$ also acts on $C^*(T^n; C(\pi^*X, T^k))$ from the right by 
\[
(\sigma \cdot a)(u_1,\ldots u_q,\tilde{x}):=\sigma (\rho (a)(u_1),\ldots ,\rho (a)(u_q), \phi_{\pi_1}(a)(\tilde{x}))
\]
for $a\in \pi_1(B_X)$, $\sigma \in C^q(T^n; C(\pi^*X, T^k))$, $u_1, \ldots , u_q\in T^n$ and $\tilde{x}\in \pi^*X$. 
\begin{lem}
The $\pi_1(B_X)$-action on $C^*(T^n; C(\pi^*X, T^k))$ commutes with the coboundary operators $\delta$. 
\end{lem}
\begin{proof}
It can be checked by using the equation $\phi_T(\rho (a)(u))\circ \phi_{\pi_1}(a)=\phi_{\pi_1}(a)\circ \phi_T(u)$ for any $a\in \pi_1(B_X)$ and $u\in T^n$. 
\end{proof}

Since $\tilde{\phi}_T(u)^{-1}\circ \tilde{\phi}_{\pi_1}(a)^{-1}\circ \tilde{\phi}_T(\rho (a)(u))\circ \tilde{\phi}_{\pi_1}(a)$ lies in $\Aut (\widetilde{P})$ for any $u\in T^n$ and $a\in \pi_1(B_X)$, the equation 
\begin{equation}\label{sigma}
\tilde{\phi}_T(u)^{-1}\circ \tilde{\phi}_{\pi_1}(a)^{-1}\circ \tilde{\phi}_T(\rho (a)(u))\circ \tilde{\phi}_{\pi_1}(a)(\tilde{p})=\tilde{p}\cdot \sigma (a, u, \pi_{\widetilde{P}}(\tilde{p}))
\end{equation}
determines the unique one-cochain $\sigma \in C^1(\pi_1(B_X); C^1(T^n; C(\pi^*X, T^k)))$. It is easy to see that $\sigma$ takes a value in $Z^1(T^n; C(\pi^*X, T^k))$. 
\begin{lem}
$\sigma$ is a cocycle. 
\end{lem}
\begin{proof}
For $a_1, a_2\in \pi_1(B_X)$, $u\in T^n$, and $\tilde{p}\in \widetilde{P}$ with $\pi_{\widetilde{P}}(\tilde{p})=\tilde{x}$, 
\[
\begin{split}
&\tilde{p}\cdot (\delta \sigma )(a_1, a_2, u, \tilde{x})\\
&=\tilde{\phi}_{\pi_1}(a_2)^{-1}\circ \tilde{\phi}_T(\rho(a_2)(u))^{-1}\circ \tilde{\phi}_{\pi_1}(a_1)^{-1}\circ \tilde{\phi}_T(\rho(a_1a_2)(u))\circ \tilde{\phi}_{\pi_1}(a_1)\circ \tilde{\phi}_{\pi_1}(a_2)\\
&\circ \tilde{\phi}_{\pi_1}(a_1a_2)^{-1}\circ \tilde{\phi}_T(\rho(a_1a_2)(u))^{-1}\circ \tilde{\phi}_{\pi_1}(a_1a_2)\circ \tilde{\phi}_T(u)\\
&\circ \tilde{\phi}_T(u)^{-1}\circ \tilde{\phi}_{\pi_1}(a_2)^{-1}\circ \tilde{\phi}_T(\rho(a_2)(u))\circ \tilde{\phi}_{\pi_1}(a_2)(\tilde{p})\\
&=\tilde{p}. 
\end{split}
\]
This proves the lemma. 
\end{proof}

Then $\sigma$ defines a cohomology class in $H^1(\pi_1(B_X);Z^1(T^n; C(\pi^*X, T^k)))$. It is easy to see that the cohomology class does not depend on the choice of liftings. We denote it by $\Fo (P)$. 
\begin{thm}\label{main-thm}
Under the assumption of Theorem~\textup{\ref{Hattori-Yoshida}}, the vanishing of $\Fo (P)$ is the necessary and sufficient condition in order that $\pi_P\colon P\to X$ admits a lifting of the local $T^n$-action $\CT$ on $X$.  
\end{thm}
\begin{proof}
If there is a lifting $\tilde{\phi}'_T\colon T^n\to \Isom_G(\widetilde{P})$ of $\phi_T$ which satisfies \eqref{semidirect-prod-condi}, then the equation 
\begin{equation}\label{tau}
\tilde{\phi}'_T(u)(\tilde{p})=\tilde{\phi}_T(u)(\tilde{p})\cdot \tau (u,\tilde{x})
\end{equation}
for $u\in T^n$ and $\tilde{p}\in \widetilde{P}$ with $\pi_{\widetilde{P}}(\tilde{p})=\tilde{x}$ defines a unique group one-cocycle $\tau \in Z^1(T^n;C(\pi^*X,T^k))$. By \eqref{sigma}, \eqref{tau}, and \eqref{semidirect-prod-condi} for $\tilde{\phi}'_T$, we can obtain the following equation 
\begin{equation}\label{coboundary}
\sigma (a, u, \tilde{x})=\tau (u,\tilde{x})\tau (\rho(a)(u), \phi_{\pi_1}(a)(\tilde{x}))^{-1}
\end{equation}
for $a\in \pi_1(B_X)$, $u\in T^n$, and $\tilde{x}\in \pi^*X$. This implies that $\Fo (P)$ vanishes.  

Conversely, suppose that $\Fo (P)$ vanishes and $\tau \in Z^1(T^n;C(\pi^*X,T^k))$ is an element which satisfies \eqref{coboundary}. Then, we define the map $\tilde{\phi}'_T\colon T^n\to \Isom_G(\widetilde{P})$ by 
\[
\tilde{\phi}'_T(u)(\tilde{p}):=\tilde{\phi}_T(u)(\tilde{p})\cdot \tau (u,\tilde{x})
\]
for $u\in T^n$ and $\tilde{p}\in \widetilde{P}$ with $\pi_{\widetilde{P}}(\tilde{p})=\tilde{x}$. It is easy to see that $\tilde{\phi}'_T$ is a homomorphism and in fact, is a lifting of the local $T^n$-action $\CT$ on $X$ to $P$. This proves the theorem. 
\end{proof}
\begin{rem}
For a general principal bundle whose structure group is not necessarily $T^k$, a lifting of a local torus action can be also defined by Definition~\ref{defn-lifting}, and Theorem~\ref{main-thm} holds without modifications not only for principal $T^k$-bundles but also for principal bundles with any abelian structure groups provided that liftings of local torus actions to the principal bundles exist. 
\end{rem}

\bibliographystyle{amsplain}
\bibliography{lifting_local}

\providecommand{\bysame}{\leavevmode\hbox to3em{\hrulefill}\thinspace}
\providecommand{\MR}{\relax\ifhmode\unskip\space\fi MR }
\providecommand{\MRhref}[2]{%
  \href{http://www.ams.org/mathscinet-getitem?mr=#1}{#2}
}
\providecommand{\href}[2]{#2}
\begin{thebibliography}{1}

\bibitem{BP}
V.~Buchstaber and T.~Panov, \emph{Torus actions and their applications in
  topology and combinatorics}, University Lecture Series, vol.~24, Amer. Math.
  Soc., Providence, RI, 2002.

\bibitem{Dav}
M.~Davis, \emph{Group generated by reflections and aspherical manifolds not
  covered by {E}uclidean space}, Ann. of Math. (2) \textbf{117} (1983), no.~2,
  293--324.

\bibitem{DJ}
M.~Davis and T.~Januszkiewicz, \emph{Convex polytopes, coxeter orbifolds and
  torus actions}, Duke Math. J. \textbf{62} (1991), no.~2, 417--451.

\bibitem{Ham}
M.~Hamilton, \emph{Quantization of toric manifolds via real polizations}, The
  talk in International Conference on Toric Topology, 2006.

\bibitem{HY}
A.~Hattori and T.~Yoshida, \emph{Lifting compact group actions in fiber
  bundles}, Japan. J. Math. (N.S.) \textbf{2} (1976), no.~1, 13--25.

\bibitem{St}
T.~E. Stewart, \emph{Lifting group actions in fiber bundles}, Ann. of Math. (2)
  \textbf{74} (1961), no.~1, 192--198.

\bibitem{Y1}
T.~Yoshida, \emph{Twisted toric structures}, arXiv:math.SG/0605376, 2006.

\bibitem{Y2}
\bysame, \emph{Local torus actions modeled on the standard represenation},
  arXiv.0710.2166, 2007.

\end{thebibliography}
\end{document}